\documentclass{amsart}

\usepackage{amsmath,amssymb,amsthm,xcolor}
\usepackage{graphicx,tikz}

\newtheorem*{thm}{Theorem}
\newtheorem*{proposition}{Proposition}

\newtheorem{lemma}{Lemma}
\newtheorem*{problem}{Problem}

\theoremstyle{definition}

\theoremstyle{remark}

\newtheorem*{ackno}{Acknowledgement}

\newcommand{\smin}{\sigma_{\min}}
\newcommand{\C}{\mathbb{C}}
\newcommand{\Z}{\mathbb{Z}}
\newcommand{\R}{\mathbb{R}}
\newcommand{\Id}{\mathrm{Id}}

\begin{document}

\title[How well--conditioned can the eigenvalue problem be?]{How well--conditioned can \\the eigenvalue problem be?}

\author[Beltrán]{Carlos Beltr\'an}
\address{Facultad de Ciencias, Universidad of Cantabria, 39005 Santander, Spain. }
\thanks{CB is Partially supported by Ministerio de
Economía y Competitividad, Gobierno de España, through grants MTM2017-83816-P,
MTM2017-90682-REDT Banco de Santander and Universidad de Cantabria
grant 21.SI01.64658.
LB is supported by the Austrian Science Fund (FWF) and the
  German Research Foundation (DFG) through the joint project
  FR 4083/3-1/I 4354.
  PG is supported by the Austrian Science Fund FWF project F5503
  (part of the Special Research Program (SFB) ``Quasi-Monte Carlo Methods:
  Theory and Applications'').
  SS is partially supported by the NSF (DMS-2123224) and the Alfred
  P. Sloan Foundation (FG-2021-14114).
}
\email{beltranc@unican.es}

\author[Betermin]{Laurent B\'etermin}
\address{Faculty of Mathematics, University of Vienna, 1090 Vienna, Austria.}
\email{laurent.betermin@univie.ac.at}

\author[Grabner]{Peter Grabner}
\address{Institute of Analysis and Number Theory, Graz University of
  Technology, Kopernikusgasse 24, 8010 Graz, Austria}
\email{peter.grabner@tugraz.at}

\author[Steinerberger]{Stefan Steinerberger}
\address{Department of Mathematics, University of Washington, Seattle, USA}
\email{steiner@uw.edu}

\begin{abstract}
  The condition number for eigenvalue computations is a well--studied
  quantity. But how small can we expect it to be? Namely, which is a perfectly
  conditioned matrix w.r.t. eigenvalue computations? In this note we answer
  this question with exact first order asymptotic.
\end{abstract}

\maketitle

\section{Introduction and Result}
For a matrix $A\in\C^{n\times n}$ and an eigenpair
$(\lambda,x)\in\C\times\C^n$, the classical Schur decomposition yields an
unitary matrix $Q$ such that
\begin{equation}\label{eq:A}
  Q^HAQ=\begin{pmatrix}\lambda&w^H\\0&B\end{pmatrix},
\end{equation}
where $\cdot\,^H$ denotes the Hermitian conjugate and $w\in\C^{n-1}$ is a
vector. Denoting by $y\in\C^n$ the corresponding left eigenvector, the
condition numbers for the eigenvalue $\lambda$ and eigenvector $x$ are given by
(see for example \cite{VanLoan})
\begin{equation*}
  \kappa_\lambda(A)=\frac{\|y\|\,\|x\|}{|y^Hx|} \quad \text{and }
  \quad \kappa_x(A)=\frac{1}{\smin(B-\lambda \cdot \Id_{(n-1) \times (n-1)})},
\end{equation*}

where $\smin$ denotes the least singular value and $\Id$ is the identity
matrix. In other words, if we allow for a $\varepsilon$--size perturbation of
$A$, the eigenpair $(\hat\lambda,\hat x$) of the perturbed matrix will satisfy
\begin{equation*}
  |\hat\lambda-\lambda|\leq \varepsilon \kappa_\lambda(A)+O(\varepsilon^2),
  \quad \angle(x,\hat x)\leq \varepsilon \kappa_x(A)+O(\varepsilon^2).
\end{equation*}
If the matrix $A \in \mathbb{C}^{n\times n}$ is diagonal with pairwise
different entries $z_1,\ldots,z_n$, we have $\kappa_{z_i}(A)=1$ for all $i$ and
the eigenvector condition number admits a simpler expression
\begin{equation*}
  \kappa_{e_i}(A)=\frac{1}{\displaystyle \min_{j\neq i}|z_i-z_j|}.
\end{equation*}
In this note we investigate the following natural question: {\em how good can
  the eigenvector conditioning of a $n\times n$ matrix be?} The answer has a
concrete application in the search for good starting points for homotopy
methods for the eigenvector problem, see \cite{armen}, but it is just such a
natural and basic question that it deserves an answer on its own right!  The
answer is not trivial. For example, if $A=\Id_{n \times n}$ we clearly have
$\kappa_x(A)=\infty$ (and it is a rather well--known fact that one can perturb
the identity matrix $\Id_{n \times n}$ with very small changes to get any desired
collection of eigenvectors).  In general, one is interested in perturbations
which are relative to the size of $A$, that is, perturbations of size
$\varepsilon\|A\|_*$ where $\|A\|_*$ is either the operator or the Frobenius
norm.
\begin{problem}\label{problem}
  Which is the optimal value for the relative--error perturbation eigenvector
  conditioning of $A\in\C^{n\times n}$, that is, which matrix minimizes the
  quantity
\begin{equation*}
\kappa_{\mathrm{max},*}(A)=\max_{x}\left(\kappa_x(A)\|A\|_*\right),
\end{equation*}
where $x$ runs over all eigenvectors of $A$ and $\|A\|_*=\|A\|_{\mathrm{Frob}}$
or $\|A\|_*=\|A\|_{\mathrm{op}}$ is the Frobenius or the operator norm?
\end{problem}
We recall the unit--side triangular (sometimes called hexagonal!) lattice in
$\C\equiv\R^2$ which is the set of points of the form
\begin{equation*}
\begin{pmatrix}
  1 & 1/2 \\
  0 & \sqrt{3}/2
\end{pmatrix}\begin{pmatrix}
               a \\
               b
             \end{pmatrix},\quad a,b\in\Z.
           \end{equation*}

\begin{center}
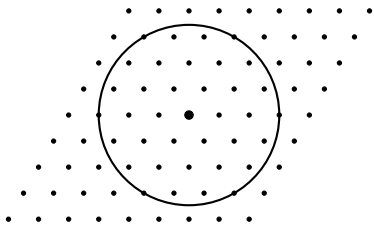
\begin{figure}[h!]
\begin{tikzpicture}[scale=0.4]
  \foreach \x in {-4,...,4}
    \foreach \y in {-4,...,4}
      {
        \filldraw (\x + 0.5*\y,0.866*\y) circle (0.07cm);
      }
  \filldraw (0,0) circle (0.14cm);
\draw [thick]   (0, 0) circle (3cm);
\end{tikzpicture}
\caption{An extremal configuration (to leading order): a circle at the origin
  and the points of a triangular lattice inside the circle.}
\label{fig:ave}
\end{figure}
\end{center}

Our main result states that the diagonal matrix whose entries range over the
triangular lattice, solves Problem \ref{problem}, at least with respect to the
leading order term.

\begin{thm}
  Let $\left\{z_1,\ldots,z_n\right\} \subset \mathbb{C}^{n}$ be first $n$
  points in the unit--side triangular lattice, in increasing modulus order (if
  two points have the same modulus, we take any of them). Then, as
  $n\to +\infty$,
\begin{align*}
  \kappa_{\mathrm{max,Frob}}(\emph{Diag}(z_1,\ldots,z_n))=&
  \frac{3^{1/4}}{2 \sqrt{\pi}} n + o(n),\\
  \kappa_{\mathrm{max,op}}(\emph{Diag}(z_1,\ldots,z_n))=&
  \frac{3^{1/4}}{\sqrt{2\pi}} \sqrt{n} + o(\sqrt{n}).
\end{align*}
Moreover, this diagonal matrix is asymptotically optimal in the sense that for
any matrix $A\in\C^{n\times n}$ the equalities above give lower bounds for
$\kappa_{\mathrm{max,Frob}}$ and $\kappa_{\mathrm{max,op}}$.
\end{thm}

It would naturally be interesting to have a better understanding of the error
terms and, in particular, to have a better understanding of the extremal
configurations. We believe it to be conceivable that our construction may
perhaps be quite close to optimal even with respect to lower order error
terms. For a configuration minimizing the Frobenius norm it is clear that the
center of mass has to be at $0$. This can of course be achieved by a
translation of less than unity and does not affect the asymptotic main term.
In order to get the Theorem, we will prove the following asymptotic inequality
which is interesting on its own right.
\begin{proposition}\label{prop:main}
  Let $p>0$ be fixed. We have, for any
    $z_1,\ldots,z_n\in \mathbb{C}$, as $n \rightarrow \infty$,
  \begin{equation*}
    \frac{1}{\displaystyle \min_{i \neq j} |z_i - z_j|}
    \left( \sum_{i=1}^{n} |z_i|^p\right)^{1/p}
    \geq \left(\frac{2}{p+2}\right)^{1/p}
    \frac{3^{1/4}}{\sqrt{2\pi}} n^{\frac12+\frac1p} +
    o\left(n^{\frac12+\frac1p}\right),
\end{equation*}
and this bound is matched by $z_1,\ldots,z_n$ as in the Theorem.
\end{proposition}

Just as in the Theorem, it might be interesting to get a better understanding
of the lower-order terms. If optimal configurations are indeed close to a
subset of the hexagonal lattice, then this is strongly related to the Gauss
Circle Problem and these techniques might apply (we observe that the function
$z \rightarrow |z|^p$ is also smoother than the cut-off function used in the
Gauss circle problem).

\section{Proofs}
\S 2.1. gives a proof of the Theorem (assuming the
Proposition). Section~\ref{sec:lemma} contains a simple geometric Lemma, the
proof of the Proposition is given in Section~\ref{sec:proof-proposition}.

\subsection{Proof of the Theorem}\label{sec:proof-theorem}
\begin{proof}
  From \eqref{eq:A} it is clear that the eigenvector conditioning of a matrix
  is invariant under conjugation by unitary matrices. From the Schur
  decomposition, we can thus assume that $A$ is upper--triangular. Now, for any
  eigenvector $x$ the definition of $\kappa_x(A)$ does not involve $w$ in
  \eqref{eq:A}, but $w$ contributes to the norm $\|A\|_*$, so the value of
  $\kappa_{\mathrm{max,*}}$ does not grow by setting $w=0$ for all
  eigenvectors. It follows that the matrix with optimal value of
  $\kappa_{\mathrm{max,*}}$ can be chosen diagonal (of course, conjugating it
  by any unitary matrix we get a normal matrix with identical
  conditioning). Thus, to prove the last claim of the Theorem we can assume
  that $A$ is diagonal, but in this case we note that
\begin{align*}
  \kappa_{\mathrm{max,Frob}}(\mathrm{Diag}(z_1,\ldots,z_n))=
  &\max_{i\neq j}\frac{\|(z_1,\ldots,z_n)\|_2}{|z_i-z_j|},\\
  \kappa_{\mathrm{max,op}}(\mathrm{Diag}(z_1,\ldots,z_n))=
  &\max_{i\neq j}\frac{\|(z_1,\ldots,z_n)\|_\infty}{|z_i-z_j|},
\end{align*}
and the result is immediate from Proposition \ref{prop:main} for the cases
$p=2$ and $p=\infty$. \end{proof}

\subsection{A Lemma}\label{sec:lemma}
We denote a ball of radius $r$ centered in the origin by
\begin{equation*}
B_r=\{z\in\C:|z|<r\}.
\end{equation*}
We say that a set $\left\{z_1, \dots, z_n\right\} \subset \C$ is $1$-separated
if, for all $i \neq j$
\begin{equation*}
| z_i - z_j| \geq 1.
\end{equation*}
We also introduce the counting function $N:[0, \infty] \rightarrow \mathbb{N}$
as follows: $N(r)$ is the cardinality of the largest possible $1$-separated set
contained in $B_r$. This quantity was already studied by L. Fejes T\'oth in
1940 who determined its growth.

\begin{lemma}[Fejes T\'oth \cite{Toth}]\label{lem}
We have, as $r\to +\infty$,
\begin{equation*}
N(r) =\frac{2 \pi}{\sqrt{3}} r^2 + o(r^2).
\end{equation*}
Equality is attained for the unit--side triangular lattice in $B_r$.
\end{lemma}

We give a simple sketch why this would be the case -- the simple geometric
argument makes use of Apollonian Circle Packings and the well-understood fact
that the asymptotically densest packing of circles in the plane is given by the
hexagonal lattice (which, not entirely coincidentally, is also a result of
Fejes T\'oth \cite{Toth2}).

\begin{proof}[Sketch of Proof] The claim on the triangular lattice (which
  implies the lower bound for the equality in the lemma) follows from
  \cite[Eq. (1.6)]{Lax}: for any lattice $L\subset \mathbb{R}^2$ (given by the
  set of points of the form $T\binom{a}{b}$ with $T\in\R^{2\times2}$ and
  $a,b\in\Z$,
\begin{equation*}
  \#\{ L\cap B_r\} = \frac{\pi r^2}{\det(T)} + O(r^{\frac{1}{3}}\sqrt{\log r}),
  \quad r\to +\infty.
\end{equation*}
In the case of the triangular lattice, $\det(T)=\sqrt{3}/2$ and we get the
desired result. As for the upper bound, we argue by contradiction. Suppose that
there exists $\varepsilon > 0$ such that there exists a divergent sequence
$(r_n) \rightarrow \infty$ such that
\begin{equation*}
N(r_n) \geq \left( \frac{2 \pi}{\sqrt{3}} + \varepsilon\right) r_n^2.
\end{equation*}
Let us pick a sufficiently large $r_n$ and let us tile $\mathbb{R}^2$ using
balls of radius $r_n$ in the fashion of a square or hexagonal lattice (it does
not really matter). This allows us to cover a large amount of space with very
efficient $1$-separated point sets.

\begin{figure}[h!]
\begin{tikzpicture}[scale = 0.3]
 \foreach \j in {0,...,8}
       {\draw (2*\j,0) circle [radius=01];}
 \foreach \j in {0,...,8}
       {\draw (2*\j,2) circle [radius=01];}
 \foreach \j in {0,...,8}
       {\draw (2*\j,4) circle [radius=01];}

 \foreach \j in {0,...,7}
       {\draw (2*\j+1,1) circle [radius=0.43];}

 \foreach \j in {0,...,7}
       {\draw (2*\j+1,3) circle [radius=0.43];}
\end{tikzpicture}
\caption{Using more disks with smaller radii allows for more precise
  approximation.}
\label{fig:pack}
\end{figure}
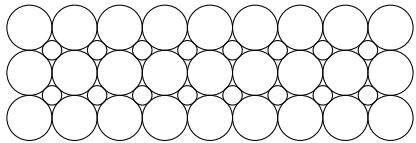

We then refine the disk packing by packing disks of smaller radius (see
Fig. \ref{fig:pack}) between the big balls and use a standard hexagonal lattice
to fill those disks. We can iterate this construction until we capture
$\sim 1 - $ of the entire plane. Not that, since $r_n$ can be chosen to be
arbitrarily large, this can always be done with a finite number of steps only
depending on $\varepsilon$. However, this then allows us to generate a period
packing of disks whose asymptotic density exceeds $\pi/\sqrt{12}$ which is a
contradiction.

\end{proof}

\subsection{Proof of the Proposition}\label{sec:proof-proposition}
\begin{proof}
  Let $p>0$ and let $\left\{z_1, \dots, z_n\right\} \subset \mathbb{R}^2$ be
  a $1$-separated set. Our goal is to prove a lower bound on
  $$ \mbox{the quantity} \qquad     \frac{1}{\displaystyle \min_{i \neq j} |z_i - z_j|}
    \left( \sum_{i=1}^{n} |z_i|^p\right)^{1/p}.$$
  We use invariance under dilation to assume without loss of generality
  that
  $$ \min_{i \neq j} |z_i - z_j| = 1$$
  and this will be assumed in all subsequent arguments.
  We abbreviate
      \begin{equation*}
M = \max_{1 \leq i \leq n} |z_i|
\end{equation*}
and write
\begin{align*}
\sum_{i=1}^{n} |z_i|^p &= p\sum_{i=1}^{n} \int_0^{|z_i |} y^{p-1}\, dy \\
&= p\int_0^M y^{p-1} \cdot \# \left\{1 \leq i \leq n: |{z_i}| > y \right\}\, dy\\
&=p\int_0^M y^{p-1} \cdot \left(n-\# \left\{1 \leq i \leq n: |{z_i}| \leq y
  \right\}\right)\, dy\\
&\geq p\int_0^My^{p-1}\max(n-N(y),0)\,dy.
      \end{align*}
      Now, by the Lemma (or the theorem given in \cite{Toth}) for every
      $\varepsilon>0$ there is an $n_0$ such that for every $n\geq n_0$
      \begin{equation}\label{eq:Mlower}
        M>\frac{3^{1/4}-\varepsilon}{\sqrt{2\pi}}\sqrt{n}.
      \end{equation}
      From this we derive
      \begin{align*}
        \sum_{i=1}^{n} |z_i|^p&\geq
        p\int_0^{\frac{3^{1/4}-\varepsilon}{\sqrt{2\pi}}\sqrt{n}}
        y^{p-1}\max\left(n-\frac{2\pi}{\sqrt3}y^2+o(y^2),0\right)\,dy\\
        &=\left(\frac{3^{1/4}-\varepsilon}{\sqrt{2\pi}}\sqrt n\right)^p
        \left(1-\frac p{p+2}\frac{(3^{1/4}-\varepsilon)^2}{\sqrt 3}\right)n+
        o\left(n^{\frac{p+2}2}\right),
      \end{align*}
which implies the proposition. For $p=\infty$, the assertion is just equation
\eqref{eq:Mlower}.

By observing that the bounds used for $N(r)$ are sharp to leading order for the
hexagonal lattice, we see that all the inequalities are asymptotically sharp.

However, there is also a direct geometric argument, valid for $p>1$, using the construction in
\cite[Lemma 2.27]{armen}: For any given $n$, we take the first $n$ points in
the unit-side triangular lattice, with increasing norm (if two points have the
same norm, we take any of them). We first note that the
function $z\to|z|^p$ is convex, hence by Jensen's inequality, for any regular
hexagon $H$ centered at $z_i$ we have
 \begin{equation*}
\mathbb{E}_{z\in H}(|z|^p)\geq |\mathbb{E}_{z\in H}(z)|^p=|z_i|^p,
\end{equation*}
where $\mathbb{E}$ has to be understood as the expected value in the
hexagon. Since by construction the points are in the unit-side triangular
lattice, the Voronoi cells surrounding each $z_i$ is a hexagon $H_i$ and we
thus have
  \begin{equation*}
    \sum_{i=1}^{n} |z_i|^p\leq \sum_{i=1}^{n} \mathbb{E}_{z\in H_i}(|z|^p)=
    \frac{1}{\mathrm{vol}(H)}\int_{\bigcup_i H_i}|z|^p\,dz.
\end{equation*}
(note that all the hexagons have the same area, which we denote by
$\mathrm{vol}(H)=\sqrt{3}/2$). Now, from the Lemma it follows that in the disk
of radius $r$ there are at least $2\pi r^2/\sqrt{3}+O(r^{1/3}\sqrt{\log r})$
points of the unit--side triangular lattice. Therefore, it follows that all the
$H_i$ are contained in a disk of radius
$r_n:=3^{1/4}\sqrt{n}/\sqrt{2\pi}+o(n)$, which yields
\begin{align*}
  \sum_{i=1}^{n} |z_i|^p&\leq \frac{1}{\mathrm{vol}(H)}
  \int_{|z|\leq r_n}|z|^p\,dz\\
  &=\frac{2\pi r_n^{p+2}}{(p+2)\mathrm{vol}(H)}=
  \frac{4\pi3^{\frac{p+2}{4}}}{\sqrt{3}(p+2)(2\pi)^{1+p/2}}n^{1+p/2}+o(n^{1+p/2}),
\end{align*}
which proves the upper bound.
\end{proof}

\begin{ackno}
  This work was initiated during the workshop ``Minimal energy problems with
  Riesz potentials'' held at the American Institute of Mathematics in May~2021.
\end{ackno}


\begin{thebibliography}{1}

\bibitem{armen}
D.~Armentano, C.~Beltr\'{a}n, P.~B\"{u}rgisser, F.~Cucker, and M.~Shub,
  \emph{{A} stable, polynomial-time algorithm for the eigenpair problem}, J.
  Eur. Math. Soc. (JEMS) \textbf{20} (2018), 1375--1437.

\bibitem{Toth}
L.~Fejes, \emph{\"{U}ber einen geometrischen {S}atz}, Math. Z. \textbf{46}
  (1940), 83--85.

\bibitem{Toth2}
L.~Fejes, \emph{\"{U}ber die dichteste {K}ugellagerung}, Math. Z. \textbf{48}
  (1943), 676--684.

\bibitem{Lax}
P.~D. Lax and R.~S. Phillips, \emph{{T}he asymptotic distribution of lattice
  points in {E}uclidean and non-{E}uclidean spaces}, J. Functional Analysis
\textbf{46} (1982), 280--350.

\bibitem{VanLoan}
C.~Van~Loan, \emph{{O}n estimating the condition of eigenvalues and
  eigenvectors}, Linear Algebra Appl. \textbf{88/89} (1987), 715--732.

\end{thebibliography}
\end{document}